\documentclass[11pt,a4paper]{article}

\usepackage{fontenc}
\usepackage[utf8]{inputenc}
\usepackage{lmodern}
\usepackage[left=2.5cm, right=2.5cm, top=2.0cm, bottom=2.5cm ]{geometry}
\usepackage{microtype}
\usepackage{amsmath,amssymb,amsthm}
\usepackage{graphicx}
\usepackage{subcaption}
\usepackage{booktabs}
\usepackage{multirow}
\usepackage{algorithm}
\usepackage{algpseudocode}
\usepackage{enumitem}
\usepackage{float}
\usepackage{tikz}
\usetikzlibrary{positioning,quotes,arrows}
\usepackage{url}
\usepackage[colorlinks=true,linkcolor=blue,citecolor=blue,urlcolor=blue]{hyperref}
\usepackage[nameinlink,noabbrev]{cleveref}
\usepackage[section]{placeins}

\newcommand{\tol}{\textsc{tol}}
\newenvironment{linenomath*}{}{}

\title{Tolerance Proportionality and Computational Stability in\\Adaptive Parallel-in-Time Runge--Kutta Methods}
\author{
Imre Fekete\textsuperscript{1,2}
\and Ferenc Izsák\textsuperscript{1,3}
\and Vendel P. Kupás\textsuperscript{1}
\and Gustaf Söderlind\textsuperscript{4}}
\date{}

\begin{document}
\maketitle

\begin{center}
\small
\textsuperscript{1} Department of Applied Analysis and Computational Mathematics, ELTE Eötvös Loránd University, Budapest, Hungary\\
\textsuperscript{2} Department of Network and Data Science, Central European University,  Vienna, Austria\\
\textsuperscript{3} NumNet HUN-REN--ELTE Research Group, Eötvös Loránd University, Budapest, Hungary\\
\textsuperscript{4} Department of Numerical Analysis, Center for Mathematical Sciences, Lund University, Lund, Sweden\\[1mm]
\end{center}

\begin{abstract}
In this paper, we investigate how adaptive time-integration strategies can be effectively combined with parallel-in-time numerical methods for solving systems of ordinary differential equations. Our focus is particularly on their influence on tolerance proportionality. We examine various grid-refinement strategies within the multigrid reduction-in-time (MGRIT) framework. Our results show that a simple adjustment to the original refinement factor can substantially improve computational stability and reliability. Through numerical experiments on standard test problems using the XBraid library, we demonstrate that parallel-in-time solutions closely match their sequential counterparts. Moreover, with the use of multiple processors, computing time can be significantly reduced.
\end{abstract}


\section{Introduction}
Adaptive numerical methods for solving ordinary differential equations (ODEs) have become a significant topic in numerical analysis. The~primary motivation for this approach lies in the challenges posed by stiff ODEs, where sudden changes in certain components may occur, and~where explicit methods are subject to severe stability restrictions. Such behavior is common, for~example, in~models of chemical reactions that involve higher-order reaction terms. The steep gradients necessitate very fine time-stepping, which can substantially slow down the~simulation.
It is important to note that realistic chemical reaction systems often involve a large number of reaction steps with widely varying reaction-rate coefficients~\cite{houston06}. Similarly, the~semi-discretization of large-scale air pollution models~\cite{jacobson05} can result in ODE systems with a vast number of~unknowns.
To address this, it is essential to identify regions where steep gradients occur, applying small time steps in their vicinity. Several algorithms have been proposed in the literature to construct efficient adaptive refinements~\cite{brenan95}.
We begin by noting that any spatial semi-discretization of partial differential equations (PDEs) yields a large, stiff system of ordinary differential equations (ODEs). As~a result, efficient numerical methods are required---even when the reaction terms are~nonstiff.
However, heuristic algorithms for selecting time steps can be prone to instabilities. For~instance, it is important to adjust the time steps smoothly by controlling the ratio between successive step sizes. This necessity naturally connects time-step selection with control theory. Indeed, some of the most effective adaptive time-stepping algorithms have been developed using control-theoretic approaches. Further details can be found in the classical works~\cite{hairernorsettwanner,Soderlind1988,Soderlind2002,Soderlind2003}. For a recent approach, we refer readers to \cite{margenov24}.
Another strategy for accelerating the numerical solution of ODEs is to use parallel-in-time algorithms. Over~the past decade, interest in these methods has been revitalized by the advent of massively parallel computing architectures. In addition to traditional applications involving PDEs~\cite{Frei2023,Steinstraesser2024,Janssens024,Zhen2024}, recent developments such as neural ordinary differential equations~\cite{NODE2018} and physics-informed neural networks~\cite{PINNs2019} have further broadened the scope of parallel-in-time methods, especially in machine learning and deep learning contexts~\cite{Tsai2024,Ibrahim2024,Parpas2025,Pamela2025}.
Parallel-in-time methods can be broadly categorized into four main groups: multiple shooting methods; domain decomposition and waveform relaxation methods; space-time multigrid methods; and direct parallel-in-time algorithms. A~comprehensive overview of these techniques is provided in the review paper~\cite{Gander2015} and the monograph~\cite{Gander2024}.
In this work, we focus on the multigrid reduction-in-time (MGRIT) method introduced in~\cite{Falgout2014}. MGRIT is closely related to the parareal algorithm developed in~\cite{Lions2001} and has since been extended in a variety of directions. These include applications to nonlinear parabolic equations~\cite{Falgout2017-Parabolic}, adaptive (variable-step) multistep methods~\cite{Falgout2019-Variable}, layer-parallel training of deep residual neural networks~\cite{Falgout2020-NN}, efficiency improvements via Richardson extrapolation~\cite{Falgout2021-Richardson}, and~hyperbolic problems such as linear advection~\cite{Falgout2023-Hyperbolic}.
Adaptive parallel-in-time methods have been explored in~\cite{Falgout2019-Variable,Falgout2021-Richardson}, with a primary focus on metrics such as time versus discretization error and processor versus wall-clock time. However, another important aspect of reliability is how the global error responds to changes in the specified tolerance. Ideally, for~problems with regular solutions, a~small change in the tolerance should result in a proportionally small change in the error. This property, referred to as tolerance proportionality, is fundamental to computational stability~\cite{Soderlind2006}.
In Section~\ref{sec:2}, we introduce the necessary concepts and tools—such as the MGRIT algorithm, basic sequential adaptive strategies, tolerance proportionality, and~computational stability—that motivate the modified refinement factor presented in the subsequent section. In~Section~\ref{sec:3}, as the main novelty of the paper, we demonstrate both tolerance proportionality and computational stability in our advanced parallel-in-time method, which builds upon the step-size selection strategy proposed in~\cite{Falgout2019-Variable}. The~proposed approach is validated through several test problems, including classical ODE systems and systems arising from the discretization of PDEs. All numerical results can be reproduced using the public GitHub repository~\cite{GitHub}. Finally, in~Section~\ref{sec:4}, we discuss the advantages, limitations, and~corresponding directions for future research.
\section{Preliminaries: Main Tools and~Requirements}\label{sec:2}
In this section, we introduce the fundamental components of our algorithm and outline the key practical requirements it must fulfill.
\subsection{The MGRIT~Algorithm}
\begin{figure}[!ht]
\centering\includegraphics[width=0.8\textwidth]{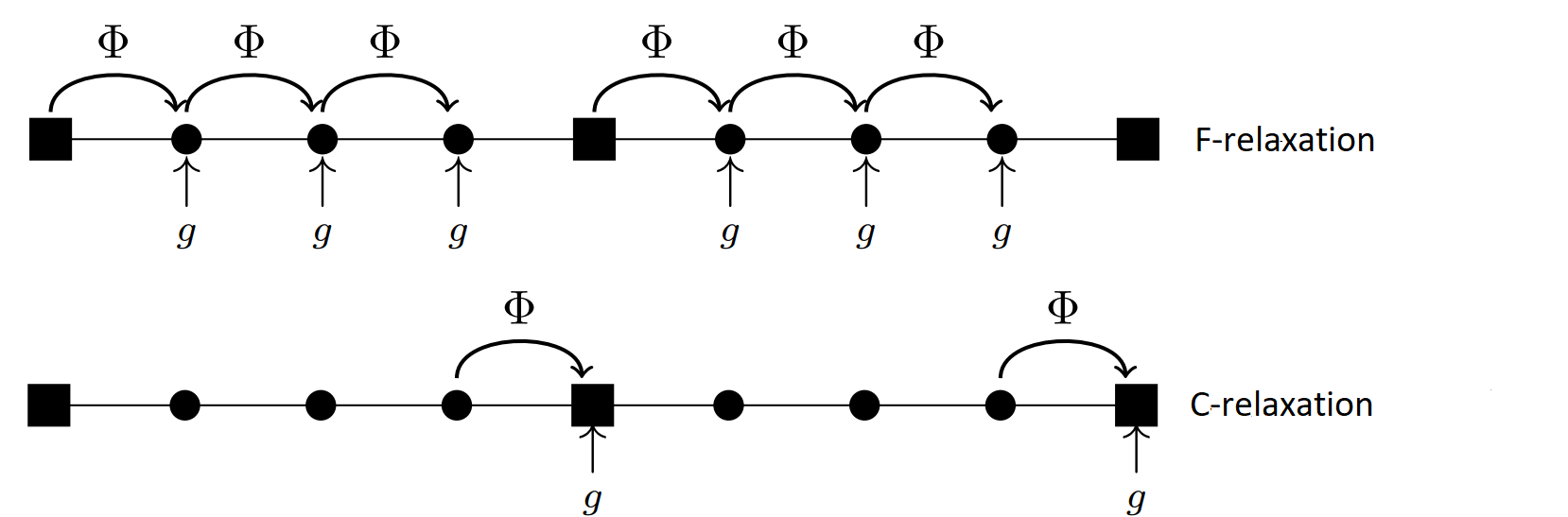}
\caption{Illustration of F- and C-relaxation~schemes.}
\label{fig:relax}
\end{figure}

\begin{figure}[!ht]
\centering\begin{subfigure}[b]{0.10\textwidth}
\centering
\includegraphics[width=\textwidth]{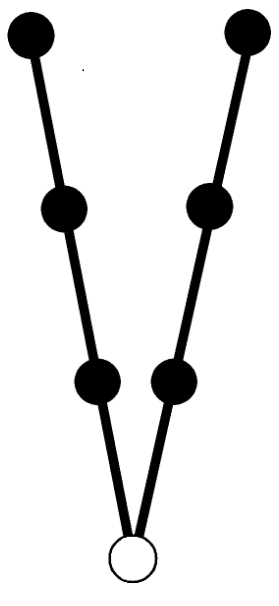}
\end{subfigure}
\begin{subfigure}[b]{0.22\textwidth}
\centering
\includegraphics[width=\textwidth]{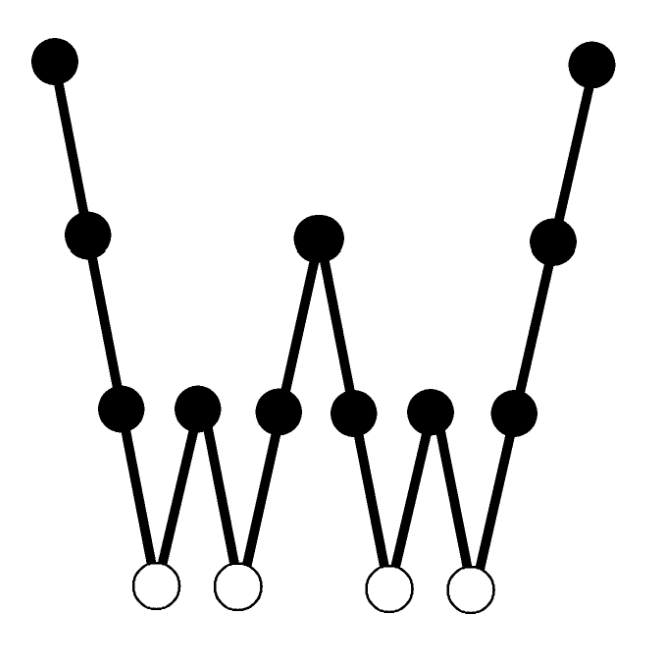}
\end{subfigure}
\begin{subfigure}[b]{0.30\textwidth}
\centering
\includegraphics[width=\textwidth]{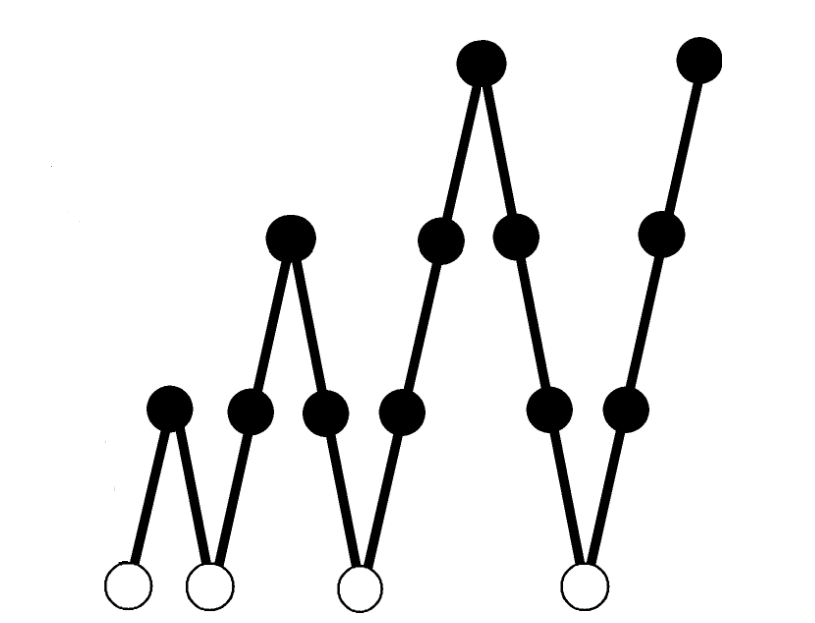}
\end{subfigure}
\caption{Multigrid cycle types: V-cycle, W-cycle, and~F-cycle.}
\label{fig:ciklusok}
\end{figure}

We provide a brief overview of the MGRIT method introduced in~\cite{Falgout2014}. Assume we are given the initial value problem
\begin{linenomath*}
\begin{align*}
u'(t) &= f(t, u(t)), \quad t \in (0, T), \\
u(0) &= u_0,
\end{align*}
\end{linenomath*}
where $f$ is a given function and $u_0$ is the initial condition. Let
$0=t_0<t_1<\ldots<t_{N-1}<t_N=T$ be a uniform mesh given on the interval $[0,T]$ with an equidistant step size $\Delta t$.
In standard time-stepping approaches, the~solution is approximated using a one-step method
\begin{linenomath*}
\begin{equation*}
\begin{split}
u_0 &= u_0, \\
u_{i+1} &= \Phi_{i+1}(u_i) + g_{i+1}, \quad i = 0, 1, \ldots, N-1,
\end{split}
\end{equation*}
\end{linenomath*}
where $\Phi$ denotes the one-step update operator, which can be practically any Runge--Kutta method, and~$g$  represents the source or correction terms. Since each time point depends on the previous one, this yields a sequential algorithm in~time.
In contrast, the~MGRIT algorithm iteratively computes an approximation to the solution over the entire time domain, starting from an initial guess. MGRIT is based on multigrid techniques and enables parallelization in the temporal domain. Assume that the right-hand side $f$ is linear, and~consider a coarser time grid $0=T_0<T_1<\ldots<T_{\frac{N}{m}-1}<T_{\frac{N}{m}}=T$ with a step size $\Delta T = m\Delta t$. On~this coarse grid, we apply another one-step method, denoted by $\Phi_{\Delta}$, which is computationally cheap and used to improve the current approximation. The~two-grid method, which is equivalent to the parareal method~\cite{Lions2001}, is outlined in Algorithm \ref{alg}. The~relaxation steps are illustrated in Figure~\ref{fig:relax}.
\begin{algorithm}[H]
    \caption{Two-grid~method.}
    \begin{algorithmic}[1]
        \State Relaxation by F-relaxation (or FCF-relaxation)
        \State Restriction: $r_{\Delta,i}=u_{mi}-(\Phi_{mi}(u_{mi-1})-g_{mi})$
        \State Solve the problem $v_{\Delta,i}=\Phi_{\Delta}(v_{\Delta,i-1})+r_{\Delta,i}$ in the coarse grid
        \State Prolongation: $u_{im} = u_{im} + v_{\Delta,i}$
        \State Relaxation by F-relaxation
    \end{algorithmic} \label{alg}
\end{algorithm}
The convergence rate of the two-grid method depends on how accurately the coarse-grid propagator $\Phi_{\Delta}$
approximates the result of performing \emph{m} successive fine-grid steps between two coarse time points.
Selecting an optimal coarse-grid propagator remains an active area of research. In~many practical settings, the~same numerical method is employed on both the fine and coarse~grids.
A potential limitation of the two-grid approach is that the coarse grid may still be too fine, rendering Step 3 of Algorithm~\ref{alg} computationally expensive. This issue can be addressed by recursively applying the same algorithm on increasingly coarser grids, following the classical multigrid paradigm. The resulting recursive structure gives rise to various cycle types, such as V-, W-, and~F-cycles, as illustrated in Figure~\ref{fig:ciklusok}. Furthermore, the algorithm can be extended to nonlinear problems using the Full Approximation Scheme~(FAS).

The MGRIT algorithm achieves optimal complexity with a runtime of $\mathcal{O}(N)$ and, unlike traditional time-stepping schemes, is massively parallelizable. In~practice, convergence to a sufficiently accurate solution is often attained within a few iterations. The numerical results of strong scaling studies on a linear problem (heat equation) and a nonlinear problem ($p$-Laplace problem) are presented in~\cite{Falgout2014,Falgout2017-Parabolic}. These results show that while the sequential one-step method outperforms MGRIT with a small number of processors, significant speedup is achieved when parallel resources are effectively~utilized.
The XBraid software package is an implementation of the MGRIT method~\cite{MGRIT}. It is written in MPI/C with interfaces for C++, Fortran 90, and~Python, and~is designed to be non-intrusive. This allows users to leverage their existing sequential time-stepping code, while XBraid manages the parallelization. Additionally, XBraid supports time refinement, enabling adaptive time-stepping strategies within the parallel-in-time framework. We used this implementation of the MGRIT method in the experiments presented in the next~section.
\subsection{Tolerance Proportionality and Computational~Stability}
If a method and its step size have been chosen properly for a given problem, it typically operates in the asymptotic regime. Thus, if~the method order is $p$, the~norm of the local error at time $t_n$ is well approximated by
\begin{linenomath*}
\begin{equation}\label{PrincipalErrorFunction}
l_n = \psi(t_n) h_n^p \, ,
\end{equation}
\end{linenomath*}
where $h_n$ is the current step size and $\psi(t)$ is the principal error function. The~latter is smooth provided that the regularity of the solution is sufficiently high. When the regularity is low, it is often preferable to use a low-order method, i.e., reduce $p$ to match the regularity of the solution so that we again have a local error of the form (\ref{PrincipalErrorFunction}) with~a smooth function $\psi(\cdot)$. In~adaptive step-size selection based on local error control, the~objective is to keep the local error constant, equal to a prescribed tolerance $\tol$, i.e.,~
\begin{linenomath*}
\begin{equation*}\label{Objective}
l_n = \psi(t_n) h_n^p = \tol \quad \Rightarrow \quad h_n = \left({\frac {\tol}{\psi(t_n)}}\right)^{1/p} .
\end{equation*}
\end{linenomath*}
It follows
 that the ``proper'' step-size sequence is smooth, and~that
\begin{linenomath*}
\begin{equation*}\label{RegularSequence}
\log h_n = {\frac {1}{p}} \big(\log \tol - \log \psi(t_n)\big).
\end{equation*}
\end{linenomath*}
Hence, if~the logarithm of the step-size sequence $\{\log h\}$ is plotted vs.\ time, a~change in $\tol$ merely shifts the graph given by $\{\log\psi\}$ vertically. In~other words, the~proper step-size sequence is a smooth function negatively correlated with the principal error function, with~$\{\log h\}$ proportional to $\log\tol$.
For these reasons, any step-size selection scheme should aim to replicate these properties. Thus, the~generated step-size sequence should be smooth, and~if well implemented, the~tolerance proportionality should be readily visible in real computations. If~the global error accumulation is regular, this implies that the global error should be proportional to the tolerance. This can, however, only be achieved with careful implementation (see, e.g.,\ \cite{Soderlind2006}). In~practice, tolerance proportionality means that the global error $u_i - u(t_i)$ is bounded below and above by the expression
\begin{linenomath*}
$$
c\cdot\tol \leq \|u_i - u(t_i)\| \leq C\cdot\tol\,,
$$
\end{linenomath*}
where the two constants satisfy $c/C \lesssim 1$. Note that this is a fairly rigorous requirement; we not only require that the local error is less than $\tol$ but~also that it is near $\tol$. Thus, if~$\tol$ is increased or reduced by a small factor, the~error must change accordingly, increasing or decreasing in proportion. This is also referred to as computational stability, as~the computational algorithm is required to produce accurate approximations in the sense that the error stays near that of previous results. In~the present context, we aim to investigate whether it is also possible to retain such properties in parallel~integration.
\subsection{Classical Recursive~Controllers}
In sequential computing, the~conventional step-size selection scheme is of the form $h_{n+1} = \rho_n\cdot h_n$, i.e.,~the step-size sequence is generated recursively along the solution trajectory. Taking logarithms, we obtain a linear recursion,
\begin{linenomath*}
\begin{equation}\label{SimpleRecursion}
\log h_{n+1} = \log h_n + \log \rho_n,
\end{equation}
\end{linenomath*}
i.e., we have a simple summation process, referred to as an integral (or I) control. The~quantity being summed is the control error $\log\rho_n$. In~the simplest case, one takes
\begin{linenomath*}
\begin{equation}\label{Elementary}
\rho_n = \left({\frac {\|e_n\|}{\tol}}\right)^{1/p},
\end{equation}
\end{linenomath*}
where $e_n$ is the local error estimate. This changes (\ref{SimpleRecursion}) into
\begin{linenomath*}
\begin{equation*}\label{SimpleRecursion1}
\log h_{n+1} = \log h_n + {\frac {1}{p}} \big( \log \|e_n\| - \log\tol \big),
\end{equation*}
\end{linenomath*}
where the deviation $\log \|e_n\| - \log\tol$ is the control error. The~step size keeps changing unless the control error is zero. Note that, in~order to generate a smooth step-size sequence, we basically advocate changing the step size every step, no matter how small the change is. Thus, older schemes of step-size halving/doubling typically generate nonsmooth step-size sequences, although~it may not be possible to avoid such sequences in connection with multigrid-type step-size~refinement.
There are many more advanced alternatives to the elementary controller (\ref{Elementary}). Some more general recursive controllers have the form
\begin{linenomath*}
\begin{equation}\label{Filter}
h_{n+1} = \rho_n^{\beta_0}\rho_{n-1}^{\beta_1}\rho_{n-2}^{\beta_2} \cdot h_n,
\end{equation}
\end{linenomath*}
where the $\rho_i$ is given by (\ref{Elementary}) and~the coefficients $\beta_0, \beta_1, \beta_2$ can be chosen to improve both stability and regularity, notably by using techniques from digital filtering (see~\cite{Soderlind2002, Soderlind2003}). Depending on the resulting structure, the~controllers may be referred to as proportional--integral (PI control) or proportional--integral--derivative (PID) control. An~important aspect is being able to reduce the gain of the controller, as~the elementary controller (\ref{Elementary}) is sometimes prone to generating oscillatory (nonsmooth) step-size sequences. There are also further possibilities that account for previous step-size sequences. The~design of the controller should be matched to the type of method and the problem class. Thus, stiff integration can benefit from special implicit methods with dedicated time-step~adaptivity.
In the present context, experiments are based on variants of (\ref{Elementary}) due to the special requirements associated with multigrid~refinement.
\section{Results}\label{sec:3}
\begin{figure}[H]
\centering\includegraphics[width=0.75\textwidth]{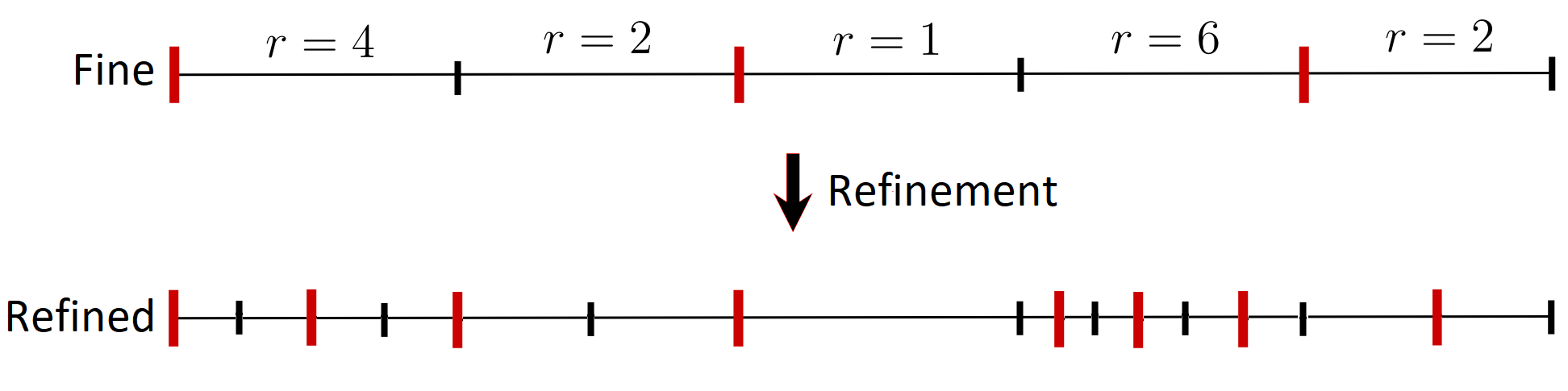}
     \caption{Schematic concept of the refinement factor $r$.}\label{fig:rfactor}
\end{figure}

The adaptivity strategy described in~\cite{Falgout2019-Variable} starts with an equidistant coarse time grid, which is further refined where necessary based on the error estimator during a given step. The~refinement factor $r$ is computed as the ratio of the current step size to the one suggested by the error estimator. Assuming that the step size--error relation is in the asymptotic regime, the~refinement factor is constructed from the formula
\begin{linenomath*}
\begin{equation}\label{eq:r_original}
    r=\frac{h_n}{h_{n+1}}= \left \lceil{ \left(\frac{\|e_n\|}{\tol}\right)^{1/p}}\right \rceil,
\end{equation}
\end{linenomath*}
where $e_n$ is the error estimate, $\tol$ is the user-provided accuracy requirement, and~$p$ is the order of the error estimator of the Runge--Kutta method. The~ceiling function is motivated by the fact that the subintervals can only be divided into an integer number of smaller subintervals. A~schematic concept of the refinement factor $r$ is shown in Figure~\ref{fig:rfactor}.
The formula \eqref{eq:r_original} without the ceiling function resembles the classical elementary I--controller \eqref{Elementary} (cf.\ \cite{hairernorsettwanner}, Chapter II.4). Therefore, we would expect tolerance proportionality properties similar to those in sequential integration, reducing the losses at the interfaces of intervals with different refinement factors. The importance of smooth step-size sequences lies in the fact that the local errors can be equidistributed, a~principle that is usually conducive to uniform global error behavior as a function of $\tol$.
\subsection{{ODE Test~Problems}}\label{subsec:ODEtestproblems}
The first test problem is the \textit{two-compartment dilution process}:
\begin{linenomath*}
\begin{equation}
\begin{aligned}
y'_1 &= -\frac{1}{5}y_1\\
y'_2 &= -\frac{2}{5}(y_2-y_1)
\end{aligned}
\label{example:linearODE}
\end{equation}
\end{linenomath*}
for $t\in[0,20]$, with~the exact solution
\begin{linenomath*}
\begin{equation*}
\begin{aligned}
y_1 &=  0.3e^{-0.2t}\\
y_2 &= 0.6(e^{-0.2t}-e^{-0.4t}).
\end{aligned}
\end{equation*}
\end{linenomath*}
The initial values  were chosen as the exact solution values at the first step. This linear problem is a standard test case with a known solution, making it ideal for assessing time-discretization accuracy and efficiency.
The second test problem is the classical nonlinear \textit{Lotka--Volterra} equation:
\begin{linenomath*}
\begin{equation}
\begin{aligned}
y'_1 &= 0.1y_1-0.3y_1y_2\\
y'_2 &= 0.5(y_1-1)y_2
\end{aligned}
\label{example:LotkaVolterra}
\end{equation}
\end{linenomath*}
for $t\in[0,62]$. The~initial values were taken as $[1,1]$.
The third problem is the classical \textit{van der Pol} equation:
\begin{linenomath*}
\begin{equation}
\begin{aligned}
y'_1 &= y_2\\
y'_2 &= \mu(1-y_1^2)y_2-y_1
\end{aligned}
\label{example:VDP}
\end{equation}
\end{linenomath*}
with $t\in[0,20]$. The~initial values were $[2,0]$. This is also a classical nonlinear test problem because it transitions from nonstiff to stiff behavior as the parameter $\mu$ increases. We chose $\mu = 2$ for the nonstiff computation.
\bigskip
\subsection{{Results for Refinement Factors $r$, $r_{mod}$, and $r_{\exp}$}}\label{subsec:oderesults}
\begin{figure}[H]
\centering
\captionsetup{skip=-5pt}
\includegraphics[width=0.7\textwidth]{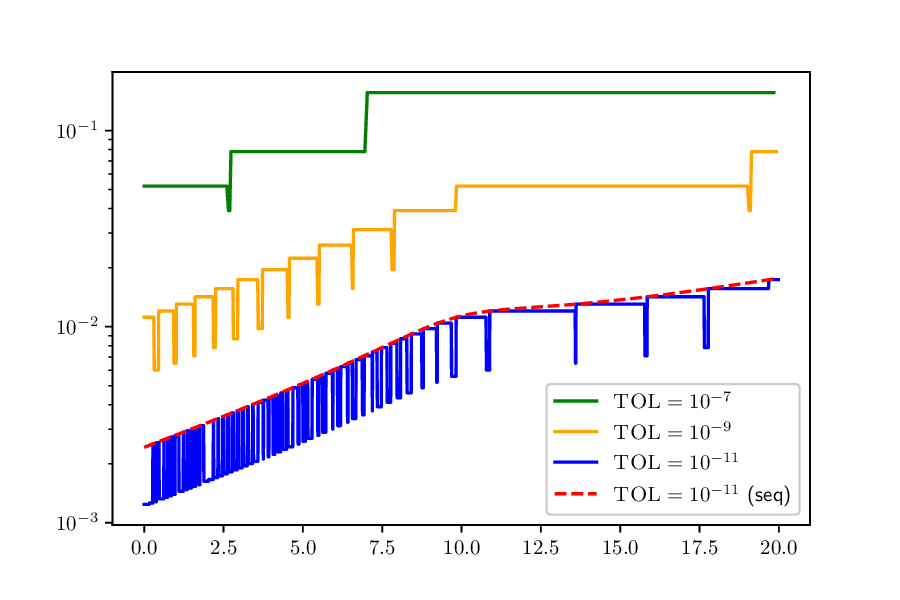}
          \caption{Step-size sequences ($y$-axis) for different TOL values in the sequential and MGRIT-based adaptive methods over the given time intervals ($x$-axis) for the linear two-compartment dilution process problem \eqref{example:linearODE}.}
        \label{fig:r_ODE_1}
\end{figure}

\begin{figure}[H]
\centering
\captionsetup{skip=-5pt}
\includegraphics[width=0.7\textwidth]{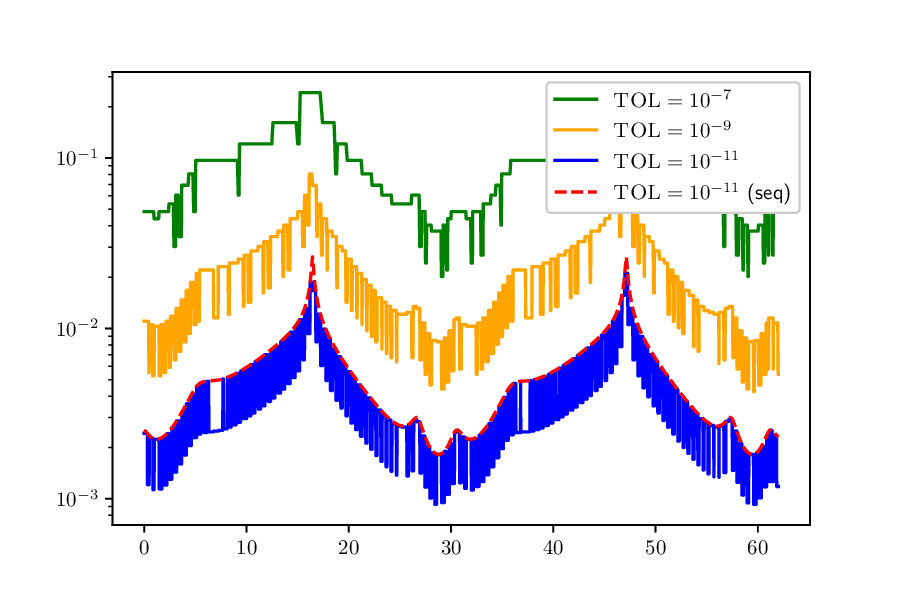}
         \caption{Step-size sequences ($y$-axis) for different TOL values in the sequential and MGRIT-based adaptive methods over the given time intervals ($x$-axis) for the Lotka--Volterra problem \eqref{example:LotkaVolterra}.}
        \label{fig:r_ODE_2}
         \end{figure}

\begin{figure}[H]
\centering
\captionsetup{skip=-5pt}
\includegraphics[width=0.7\textwidth]{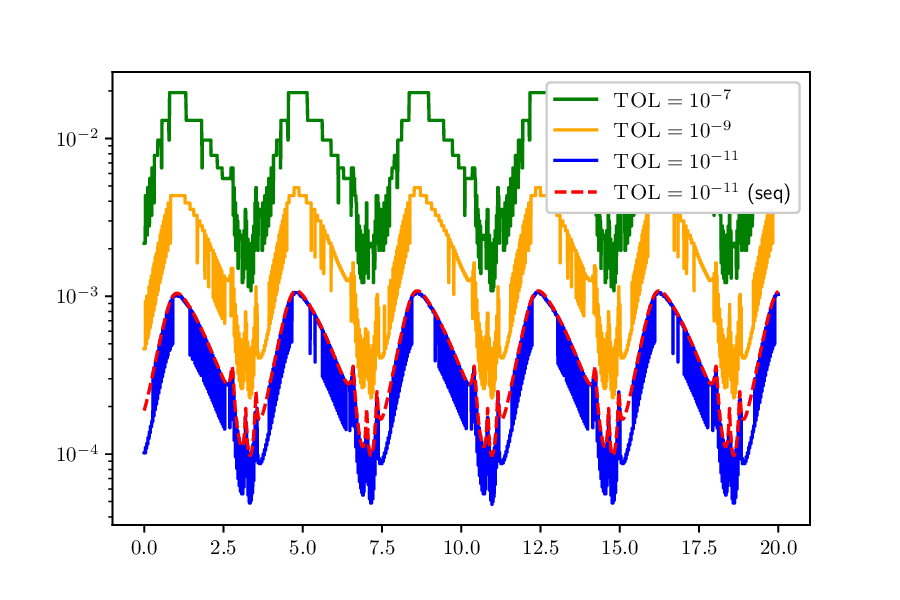}
         \caption{Step-size sequences ($y$-axis) for different TOL values in the sequential and MGRIT-based adaptive methods over the given time intervals ($x$-axis) for the van der Pol problem \eqref{example:VDP}.}
        \label{fig:r_ODE_3}
         \end{figure}

{For the ODE test problems in Section} \ref{subsec:ODEtestproblems}, smooth step-size sequences produced tolerance proportionality in the sequential case.
The MGRIT procedure, together with \eqref{eq:r_original} from the open-source XBraid library~\cite{MGRIT}, was used with the  Bogacki--Shampine method to take into account its stability~region.
The tolerance proportionality results can be seen in Figures~\ref{fig:r_ODE_1}--\ref{fig:r_ODE_3}. In~Figures~\ref{fig:r_ODE_1}--\ref{fig:r_ODE_3}, we plot the step-size sequence for each test problem. A~reasonable error tolerance proportionality is observed, as~the step size on average scales in accordance with the specified error tolerance. At~the finest tolerance, $\tol=10^{-11}$, the~MGRIT-based step-size sequence is also close to that of the sequential integration for the same~tolerance.
However, in~every test case, a~zigzag pattern appears in the step-size sequence, which contrasts with the smooth step-size sequence expected in classical sequential implementation theory~\cite{Soderlind2003}. This effect arises since the refinement factor is highly sensitive to small values due to the ceiling function. Pointwise tolerance proportionality therefore suffers, and~it is unlikely that the overall method is as efficient as it could be, since local errors are far from~equidistributed.
\vspace{-30 pt}
\vspace{-30 pt}
\FloatBarrier
We therefore conclude that it is worth considering alternatives to \eqref{eq:r_original}. Thus, we define a modified refinement factor $r_{mod}$ as follows:
\begin{equation}\label{eq:r_modified}
    r_{mod}=\frac{h_n}{h_{n+1}}= \left \lceil{ \left(\frac{\|e_n\|}{\tol}\right)^{1/p}-c}\right \rceil.
\end{equation}
Here, the~compensating  constant is taken as $c\in [0.2,0.4]$ in \eqref{eq:r_modified}.
\subsubsection{{{The Final {Algorithm} with $r_{\textrm{mod}}$}}}
\begin{figure}[H]
\centering
\vspace{30 pt}
\includegraphics[width=0.7\textwidth]{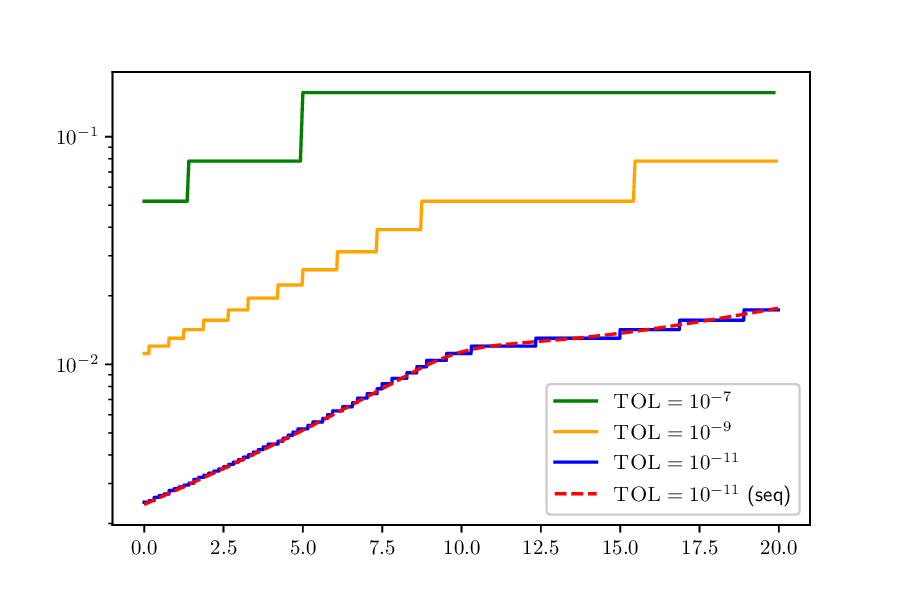}
         \caption{The experiment shown in Figure~\ref{fig:r_ODE_1} was re-simulated using the modified $r_{mod}$ refinement factor with a compensating constant $c = 0.4$.}
         \label{fig:r_mod_ODE_1}
\end{figure}

\begin{figure}[H]
\centering
\includegraphics[width=0.7\textwidth]{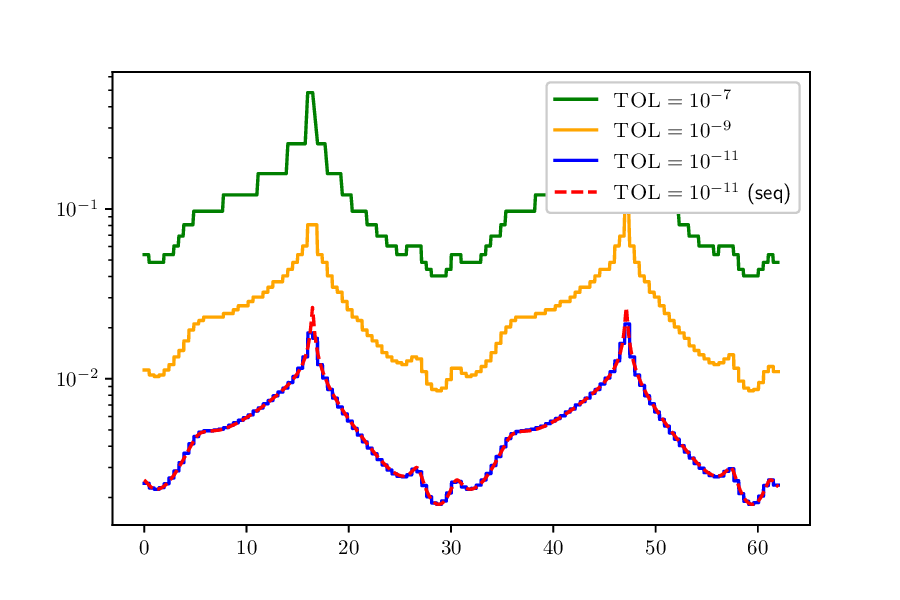}
\caption{The experiment shown in Figure~\ref{fig:r_ODE_2} was re-simulated using the modified $r_{mod}$ refinement factor with a compensating constant $c = 0.4$.}
         \label{fig:r_mod_ODE_2}
\end{figure}

\begin{figure}[H]
\centering
\vspace{-30pt}
\includegraphics[width=0.7\textwidth]{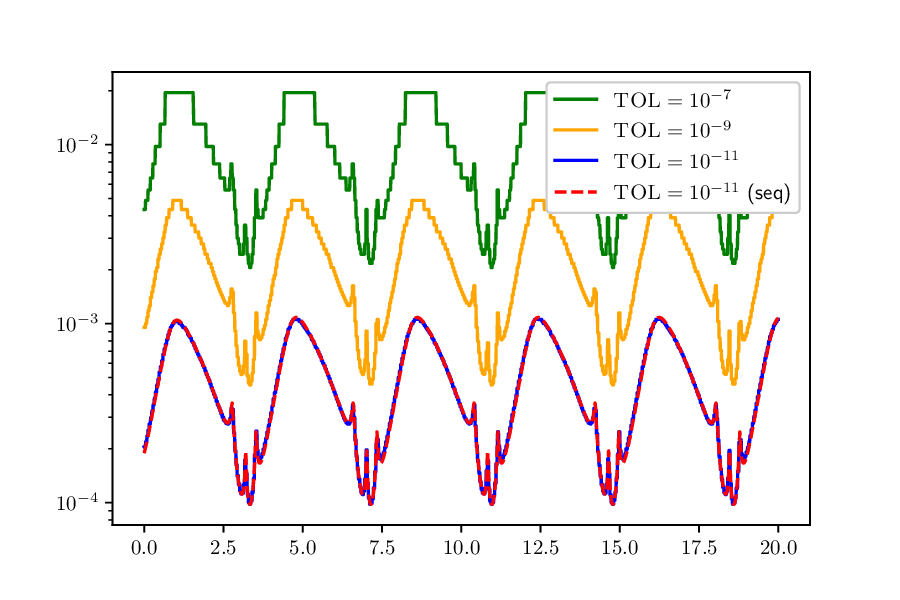}
\caption{The experiment shown in Figure~\ref{fig:r_ODE_3} was re-simulated using the modified $r_{mod}$ refinement factor with a compensating constant $c = 0.4$.}
         \label{fig:r_mod_ODE_3}
\end{figure}

Completing the adaptive MGRIT procedure with this essential modification and summarizing all steps, we apply the following algorithm:
\begin{itemize}
    \itemsep 0em
    \item [(A1)] Set an initial temporal grid and an initial guess of the solution (see Figure \ref{fig:rfactor}).
    \item [(A2)] Apply a given Runge--Kutta method in parallel from these initial values, and use multigrid cycles to improve the current approximation of the solution (see Figure \ref{fig:ciklusok}).
    \item [(A3)] Determine the refinement factor for each subinterval on the finest grid based on \eqref{eq:r_modified}.
    \item [(A4)] \textbf{{IF}}
 {no refinement occurred and the accuracy criterion is satisfied} \textbf{{then STOP}}.
    \item [(A5)] \textbf{{ELSE}} go to Step 2.
\end{itemize}
{For the tests,}
 we used the embedded Runge--Kutta methods (described later) for time stepping and error estimation. We started on a uniform temporal grid and applied three V-cycles between two consecutive refinements.
As shown in Figures~\ref{fig:r_mod_ODE_1}--\ref{fig:r_mod_ODE_3}, the~modified refinement largely suppressed the zigzag patterns, resulting in smoother step-size sequences and improved tolerance~proportionality.
\FloatBarrier
\subsubsection{{An Attempt with Exponential-Forgetting~Filters}}
An alternative idea is to reduce the integral gain in (\ref{eq:r_original}) from $1$ to, say, $2/3$, i.e., change the exponent from $1/p$ to $2/(3p)$ by constructing another modified refinement factor. As~this corresponds to an exponential-forgetting filter, we denote this new refinement factor by $r_{\exp}$. However, this does not result in a smoother step-size sequence, although~it is well known to have such an effect in the sequential case. Again, this appears to be due to the ceiling function. The corresponding figures can be found in the GitHub repository \cite{GitHub}.
\subsubsection{{Tolerance Proportionality for the Algorithm (A1)--(A5)}}
\begin{figure}[H]
\centering
\captionsetup{skip=0pt}
\includegraphics[width=0.7\textwidth]{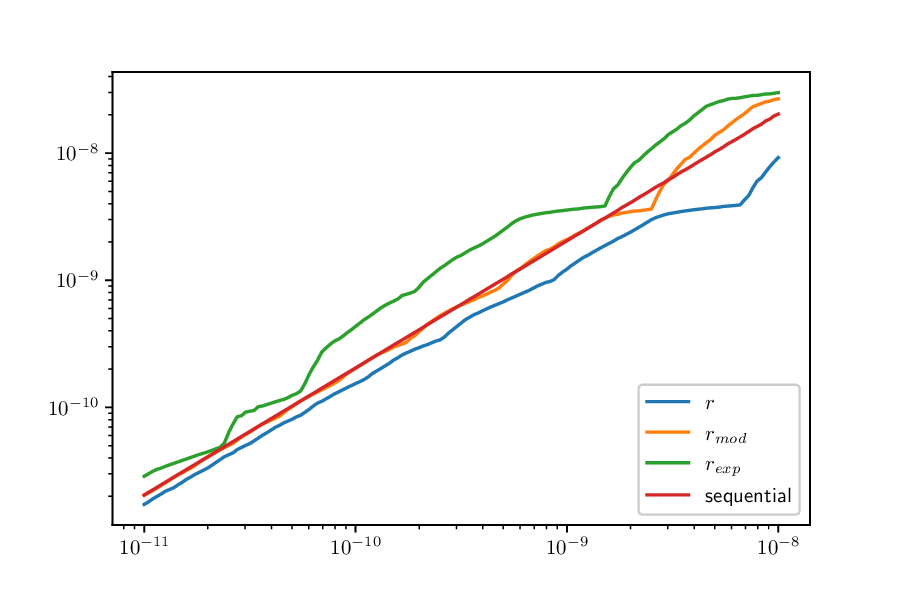}
 \caption{Achieved accuracy ($y$-axis) versus tolerance ($x$-axis) for the linear two-compartment dilution process problem \eqref{example:linearODE}, using the three refinement-factor~strategies.}
        \label{fig:accuracy_tol_ODE_1}
\end{figure}

\begin{figure}[H]
\centering
\captionsetup{skip=0pt}
\includegraphics[width=0.7\textwidth]{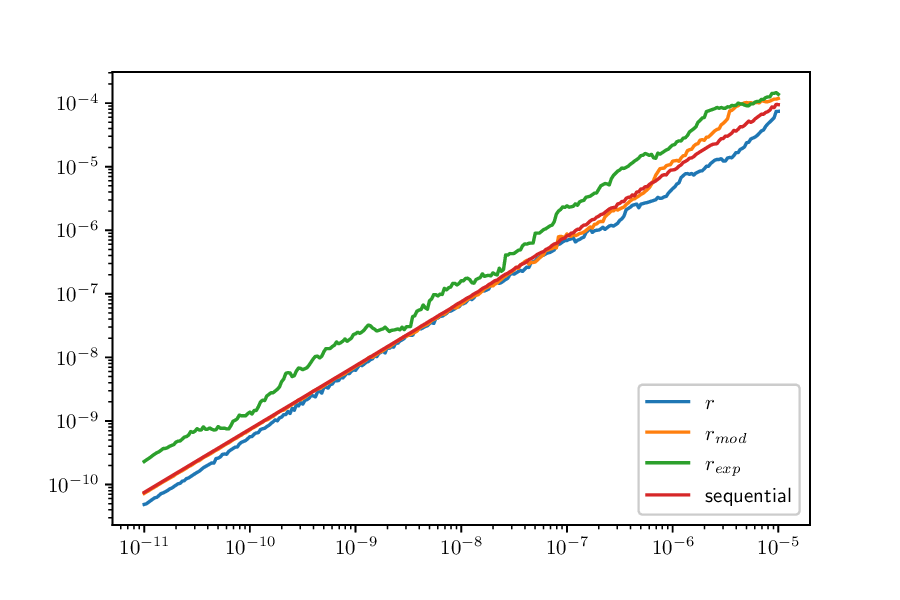}
          \caption{Achieved accuracy ($y$-axis) versus tolerance ($x$-axis) for the Lotka--Volterra problem \eqref{example:LotkaVolterra}, using the three refinement-factor~strategies.}
        \label{fig:accuracy_tol_ODE_2}
\end{figure}

\begin{figure}[H]
\centering
\captionsetup{skip=0pt}
\includegraphics[width=0.7\textwidth]{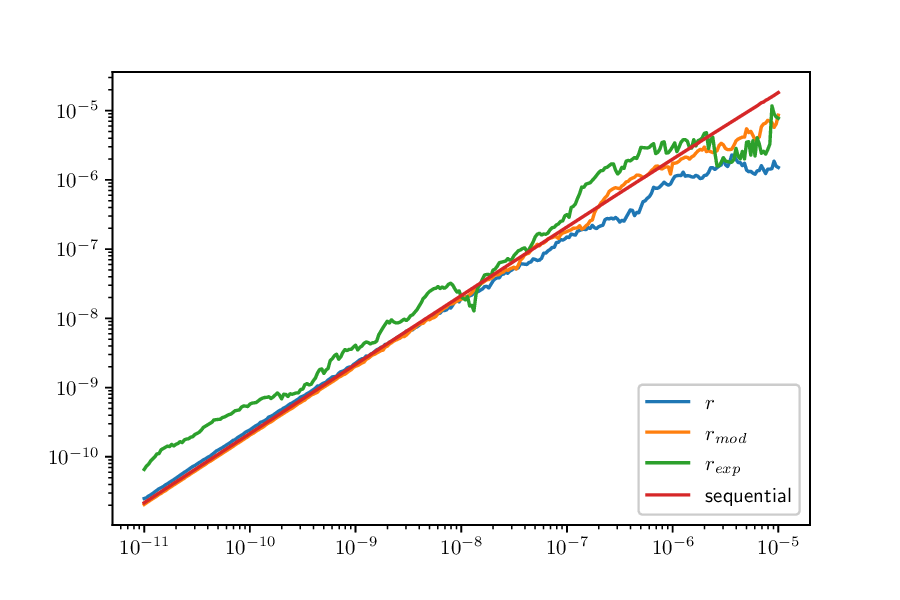}
 \caption{Achieved accuracy ($y$-axis) versus tolerance ($x$-axis) for the van der Pol problem \eqref{example:VDP}, using the three refinement-factor~strategies.}
        \label{fig:accuracy_tol_ODE_3}
\end{figure}

In connection with the concept of computational stability discussed in~\cite{Soderlind2006}, the~achieved accuracy is plotted versus $\tol$. Ideally, a~log--log plot would be expected to show a nearly straight line if the step-size sequence accurately tracks the local errors. This tolerance proportionality ensures that a reduction in $\tol$ leads to a corresponding improvement in the achieved accuracy. This accuracy versus tolerance for the three different approaches is shown in Figures~\ref{fig:accuracy_tol_ODE_1}--\ref{fig:accuracy_tol_ODE_3}.
\vspace{-30pt}
\vspace{-30pt}
Based on Figures~\ref{fig:accuracy_tol_ODE_1}--\ref{fig:accuracy_tol_ODE_3}, the~refinement factor $r_{mod}$ demonstrates the best computational stability while maintaining performance comparable to sequential computation. For~sharp tolerances, the~results are effectively identical. In~contrast, the~exponential-forgetting alternative $r_{\exp}$ provides less~benefit.
Based on Figures \ref{fig:r_ODE_1}--\ref{fig:accuracy_tol_ODE_3}, it is qualitatively clear that $r_{mod}$ is the best among the three discussed refinement factors. Naturally, one may ask whether the suppression of zigzag patterns is also beneficial in a quantitative sense. To~investigate this, we compare the $r_{mod}$ refinement factor with the less zigzag-prone $r$ refinement factor rather than with the more erratic $r_{\exp}$. As~shown in Table~\ref{tab:steps}, $r_{mod}$ makes the step-size sequences smoother, improves tolerance proportionality, and significantly reduces computational cost.
\begin{table}[H]
\centering\caption{{{Number of time}
 steps after the final refinement for different $\tol$ values, using refinement factors $r$ and $r_{mod}$, for the ODE test problems} \eqref{example:linearODE}--\eqref{example:VDP}.}
    \tabcolsep=0.49cm
\begin{tabular}{ccccccc}
\toprule
\textbf{Problem} & \textbf{\tol} & \boldmath{$10^{-5}$} & \boldmath{$10^{-7}$} & \boldmath{$10^{-9}$} & \boldmath{$10^{-11}$} & \boldmath{$10^{-13}$}\\
\midrule
\multirow{2}{*}{\eqref{example:linearODE}} & $r$ & $128$  & $191$ & $718$ & $3997$ &  $22,562$\\
                        & $r_{mod}$ & $128$ & $169$ & $620$ & $2828$ & $13,076$\\
\multirow{2}{*}{\eqref{example:LotkaVolterra}}     & $r$ & $218$ & $855$ & $4237$ & $22,371$ &  $107,982$\\
                        & $r_{mod}$ & $184$  & $754$ & $3449$ & $15,963$ & $74,047$\\
\multirow{2}{*}{\eqref{example:VDP}}    & $r$ & $931$ & $3777$ & $20,432$ & $100,043$ & $475,284$\\
                        & $r_{mod}$ & $773$ & $3070$ & $14,090$ & $65,258$ & $302,692$\\
\bottomrule
\end{tabular}
    \label{tab:steps}
\end{table}
\unskip
\subsection{{PDE Test Problems and Related~Results}}\label{subsec:pderesults}
\begin{figure}[!ht]
\centering\begin{subfigure}[b]{0.495\textwidth}
        \centering
        \includegraphics[width=0.94\textwidth, trim=10 10 10 10, clip]{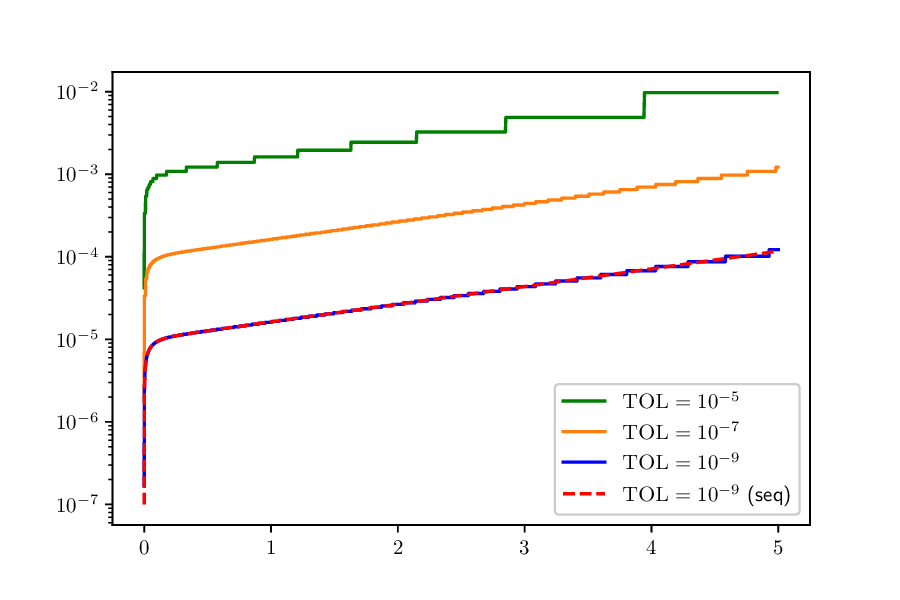}
    \end{subfigure}
    \begin{subfigure}[b]{0.495\textwidth}
        \centering
        \includegraphics[width=0.94\textwidth, trim=10 10 10 10, clip]{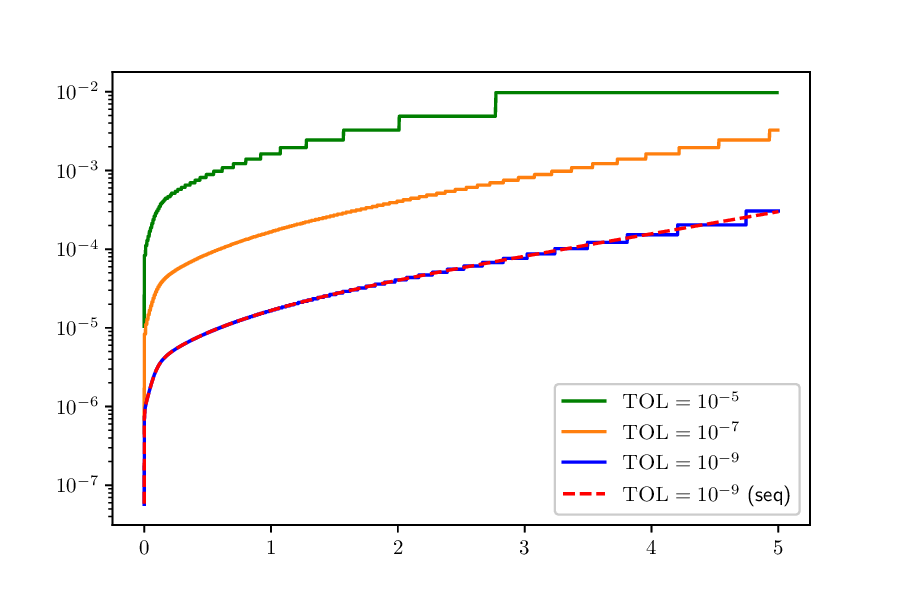}
    \end{subfigure}
    \begin{subfigure}[b]{0.495\textwidth}
        \centering
        \includegraphics[width=0.94\textwidth, trim=10 10 10 10, clip]{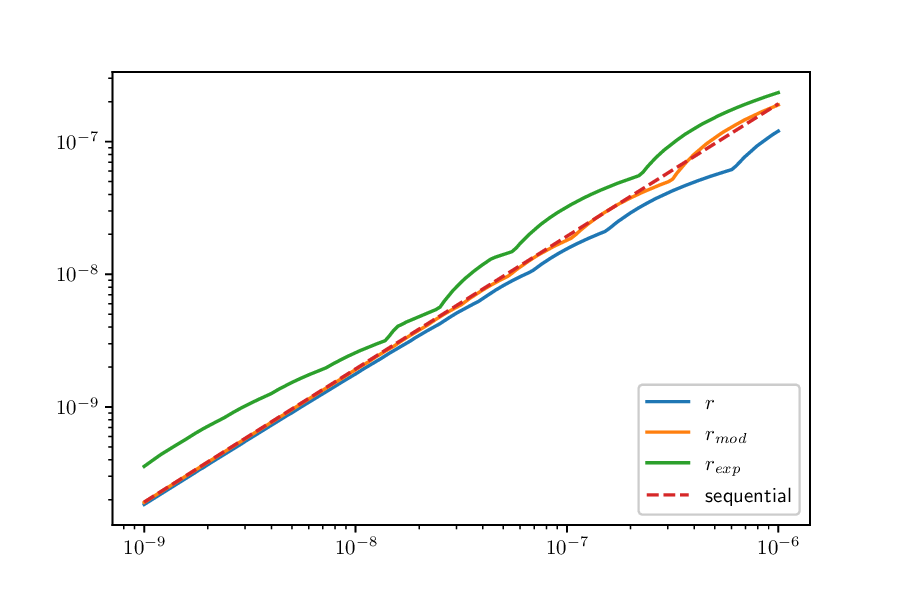}
    \end{subfigure}
    \begin{subfigure}[b]{0.495\textwidth}
        \centering
        \includegraphics[width=0.94\textwidth, trim=10 10 10 10, clip]{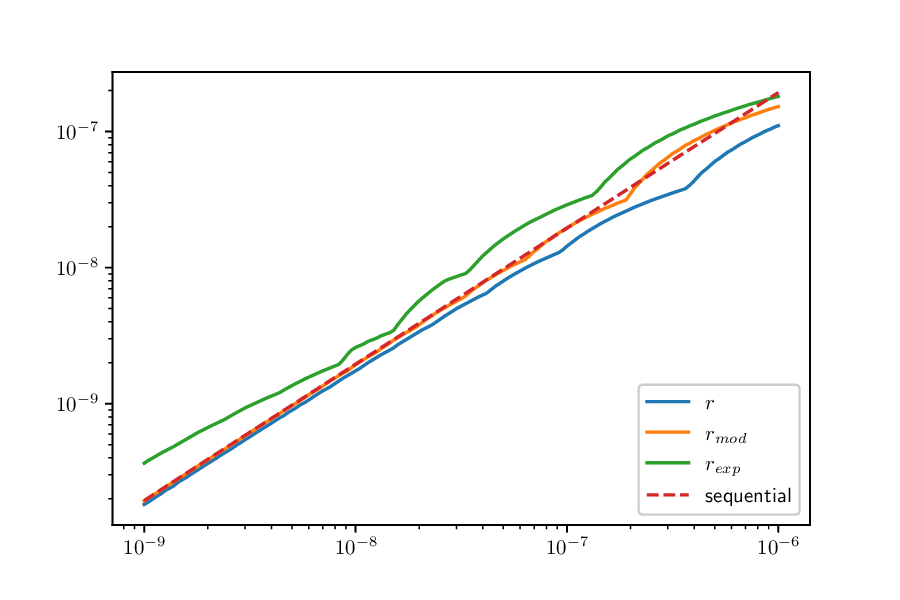}
    \end{subfigure}
    \begin{subfigure}[b]{0.495\textwidth}
        \centering
        \includegraphics[width=0.94\textwidth, trim=10 10 10 10, clip]{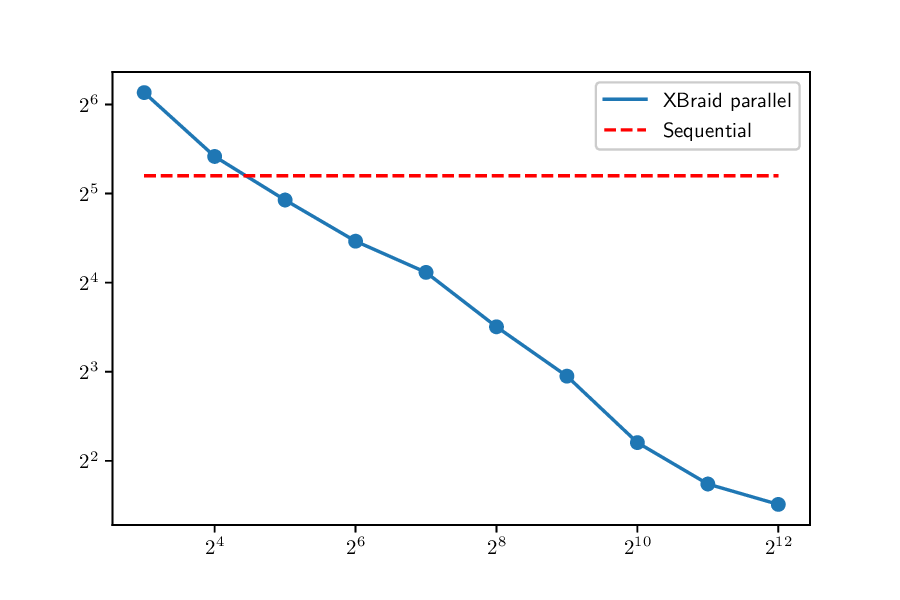}
    \end{subfigure}
    \begin{subfigure}[b]{0.495\textwidth}
        \centering
        \includegraphics[width=0.94\textwidth, trim=10 10 10 10, clip]{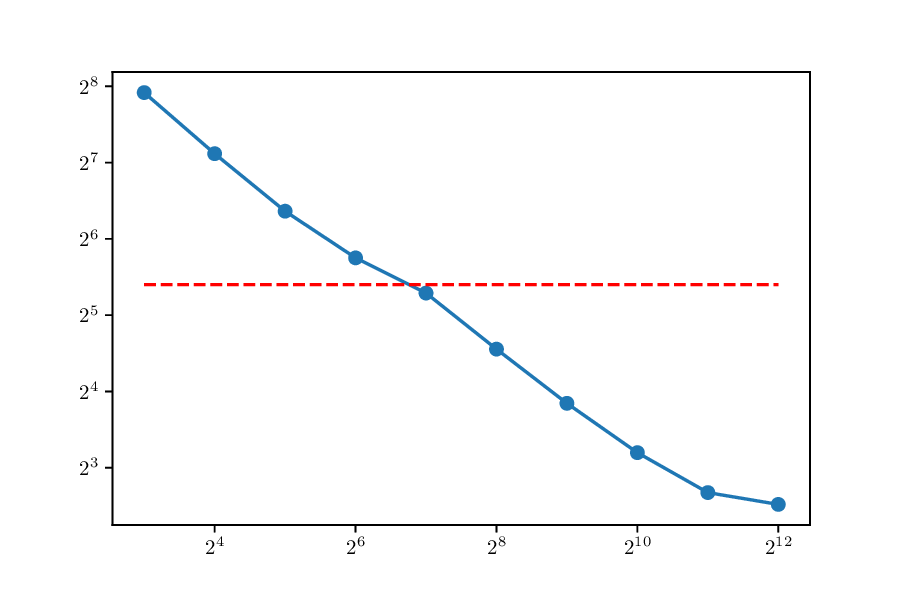}
    \end{subfigure}
    \caption{Results for the linear (left column) and nonlinear (right column) diffusion Equation~\eqref{eq:heat}. Upper panel: step-size sequences ($y$-axis) over the given time interval ($x$-axis). Center panel: achieved accuracy ($y$-axis) versus tolerance ($x$-axis). Lower panel: runtime in seconds ($y$-axis) versus number of processors ($x$-axis).}
    \label{fig:heat}
\end{figure}

For the three ODE test problems in Section \ref{subsec:ODEtestproblems}, there was no significant difference in runtime between the sequential and MGRIT-based parallel-in-time algorithms using eight processors. However, a~notable runtime advantage was observed in the case of partial differential equations (PDEs).
Our test problems for PDEs are special cases of one-dimensional, possibly nonlinear diffusion equations with Dirichlet boundary conditions. Specifically, we consider the following formulation:
\begin{equation}
    \begin{split}
        \partial_t u(t,x)&=\frac{1}{10} \partial_{xx}(u^p(t,x)), \quad x\in (0,1),\ t\in(0,5)\\
        u(0,x)&=g(x),\quad x\in (0,1)\\
        u(t,0)&=g(0)\ \text{and}\ u(t,1)=g(1),\quad t\in [0,5],\\
    \end{split}\label{eq:heat}
\end{equation}
where $g(x)=10e^{-(x-0.435576)^2}-7.27185$. The~linear case corresponds to $p=1$, while we take $p=2$ in the nonlinear case. Based on the ODE test results in Section \ref{subsec:oderesults}, we only consider the modified refinement factor $r_{mod}$ using $c=0.4$. In~space, we use the classical second-order central difference scheme. For~time integration, we use a two-stage, stiffly accurate, L-stable, singly diagonally implicit Runge--Kutta method (see~\cite{kennedycarpenter}, Butcher tableau (221), with the choices $\gamma=(2-\sqrt{2})/2$ and $\hat{b}_2=0$). When selecting the spatial and temporal step sizes, we make sure that the necessary CFL condition is satisfied, both for $p=1$ and $p=2$ in \eqref{eq:heat}. The~results related to the tolerance proportionality, achieved accuracy, and runtime can be seen in Figure~\ref{fig:heat}.

Figure~\ref{fig:heat} demonstrates that both the linear and nonlinear diffusion Equation~\eqref{eq:heat} show good tolerance proportionality and computational stability, as~well as the expected qualitative behavior of achieved accuracy versus tolerance. In~addition, the~expected reduction in runtime due to parallelization is also~observed.
\section{Discussion}\label{sec:4}
We summarize the main advantages of our proposed method and, for~the sake of completeness, also highlight the limitations of the framework, which suggest interesting research directions in their own right.
\subsection{{Advantages}}\label{subsec:advantages}
\enlargethispage{\baselineskip}
The computational results in Section \ref{sec:3} indicate that for MGRIT-based adaptive parallel-in-time Runge--Kutta methods, the~relevant quality indicators necessary for computational stability are satisfied for the refinement factors $r_{mod}$ and $r$. As we have seen in~the case of the PDE test problems in Section \ref{subsec:pderesults}, these methods significantly outperform sequential methods in terms of computational cost beyond a certain number of processors. Furthermore, based on the qualitative results of the different ODE and PDE test problems, as well as the number of time steps in Table \ref{tab:steps}, the~refinement factor $r_{mod}$ is the preferred approach both qualitatively and quantitatively. Therefore, the approach involving the refinement factor $r_{mod}$ provides a robust and reliable framework for users utilizing high-performance computing~simulations.
\subsection{Limitations}
\begin{figure}[H]
\centering\includegraphics[width=0.6\textwidth]{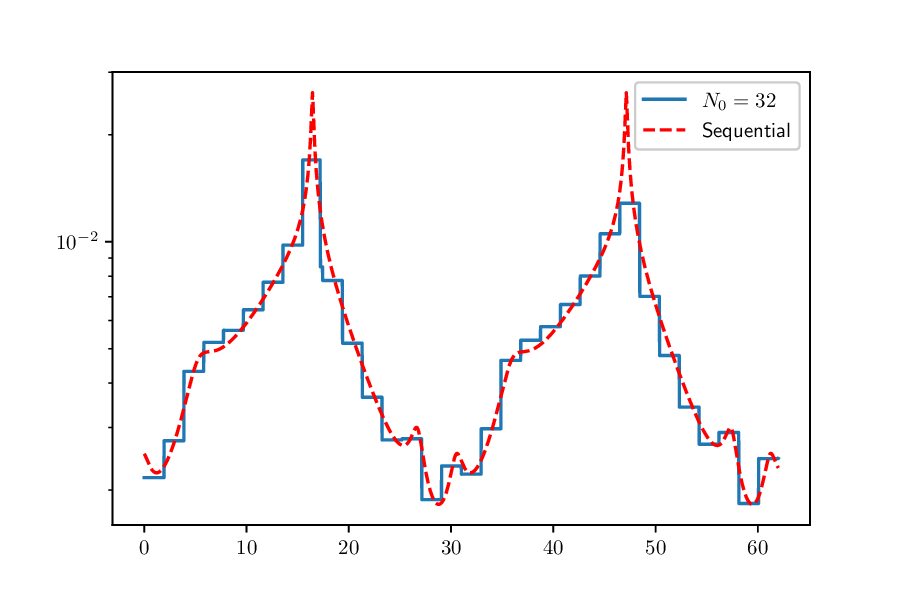}\hspace{-4mm}
    \includegraphics[width=0.6\textwidth]{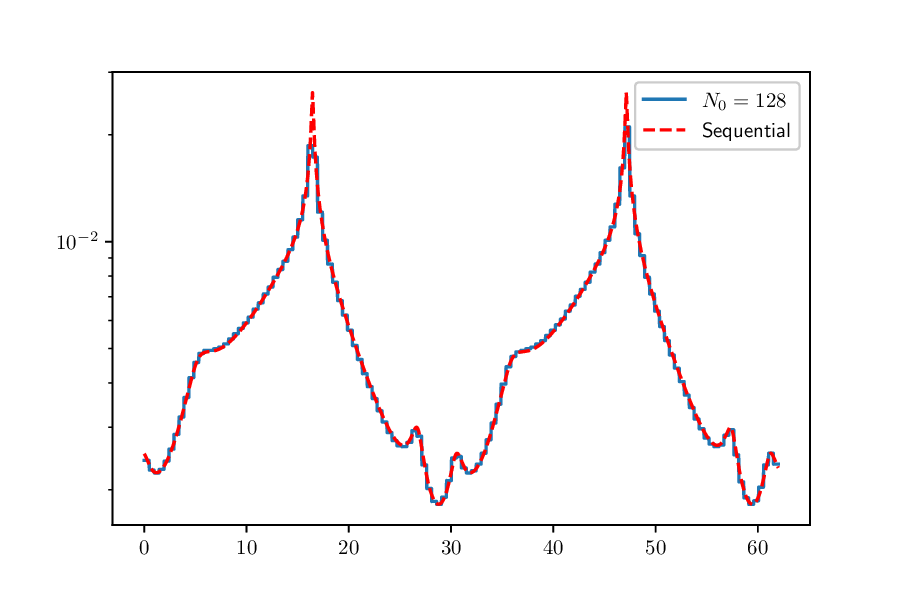}\hspace{-4mm}
    \includegraphics[width=0.6\textwidth]{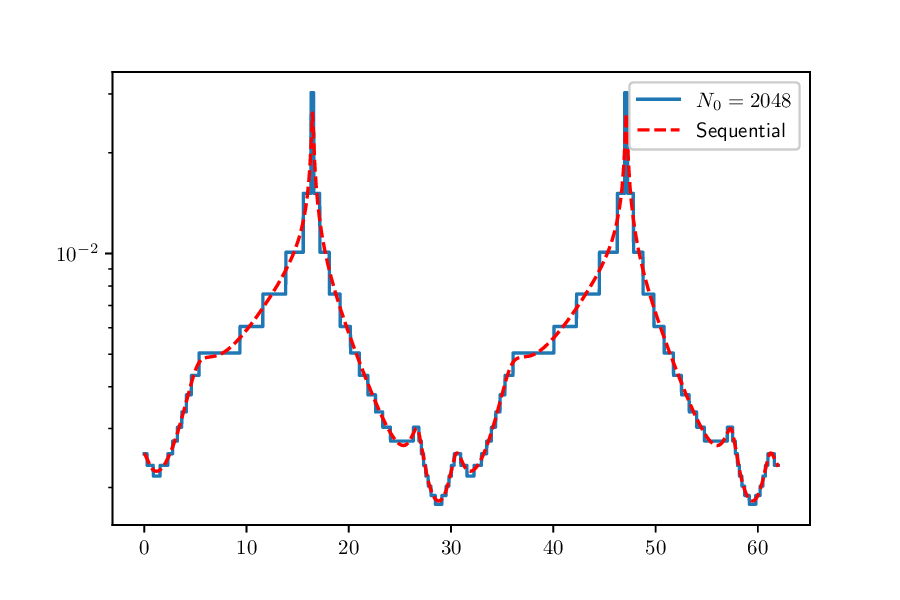}
    \caption{Effect of the initial subinterval choice on the step-size sequence ($y$-axis) over the time interval ($x$-axis) for the Lotka--Volterra Equation~\eqref{example:LotkaVolterra}.}
    \label{fig:intial_subintervals}
\end{figure}

\begin{figure}[!ht]
\centering\includegraphics[height=0.5\textwidth]{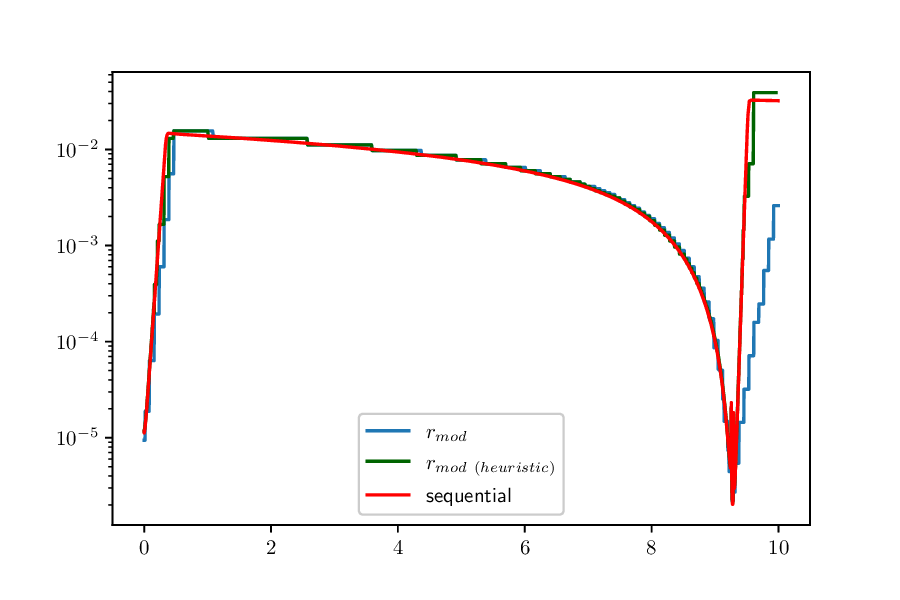}
     \caption{Step-size sequences ($y$-axis) for the stiff van der Pol Equation~\eqref{example:VDP} over the given time interval ($x$-axis). $r_{mod\ (heuristics)}$ denotes our heuristic approach for a stiff~problem.}
    \label{fig:stiff}
\end{figure}

The first limitation is that, for~each test problem, the~MGRIT algorithm must be run with an appropriate initial subinterval partitioning to achieve computational stability. Figure \ref{fig:intial_subintervals} illustrates this phenomenon. This characteristic can be seen as analogous to the choice of the initial step size in sequential algorithms (see~\cite{hairernorsettwanner}, Chapter II.4, and~\cite{Soderlind2017}, Section~4.2
).
The second limitation may be that our current results have provided satisfactory quality only for nonstiff (or at most mildly stiff) problems. This naturally suggests a research direction focusing on stiff problems. Due to their much larger variation in step sizes, they pose a more significant challenge. Based on the structure of MGRIT methods, it is clear that the analogy of PI- and PID-type controllers \eqref{Filter}, which are commonly used in sequential algorithms (for details, see~\cite{Soderlind2003}), needs reconsideration in this setting. However, from~this theoretical background, we also know that for stiff problems, PI controllers of a different parametrization are needed. An~example of this can be seen in Figure~\ref{fig:stiff}, which illustrates the case of the stiff van der Pol equation \eqref{example:VDP} with $\mu=30$ and initial value $[1.5, 3]$ on the time interval $[0,10]$. Using heuristics, we were able to significantly improve the step-size~sequence.
\vspace{-12pt}
\unskip
\subsection{{Future~Research}}
A closer examination of the above-mentioned heuristic is needed to fully understand its role in connection with parallel-in-time methods. Based on our observations, we also conclude that by incorporating a stiffness indicator (for details, see~\cite{Soderlind2015,Soderlind2024}, Chapter 21.3) as a measuring and switching mechanism, we may obtain better results for stiff problems in the current framework. Further research in this direction will lead to dedicated methods specifically designed for stiff~problems.
All simulation figures (Figures \ref{fig:r_ODE_1}--\ref{fig:stiff}) in this work are fully reproducible using the files available in the GitHub repository \cite{GitHub}.

\section*{Funding and Acknowledgements}
This research was supported by the Ministry of Innovation and Technology NRDI Office within the framework of the Artificial Intelligence National Laboratory Program RRF-2.3.1-21-2022-00004, and~within the framework of the Thematic Excellence Program ELTE TKP 2021-NKTA-62. The~authors thank Sidafa Conde (former Senior Member of Technical Staff at Sandia National Laboratories, USA) for his assistance with the initial XBraid runs. We acknowledge the Digital Government Development and Project Management Ltd. for awarding us access to the Komondor HPC facility based in Hungary.

\end{document}